\documentclass[11pt,a4paper]{article}

\usepackage{verbatim,amssymb,amsbsy,amscd,amsmath,amsthm,amsfonts,mathrsfs,bm,graphicx,fullpage}

\usepackage[mathscr]{eucal}

\newcommand{\dist}{{\rm dist}}

\newtheorem{dfn}{Definition}[section]
\newtheorem{thm}[dfn]{Theorem}

\newtheorem{rem}[dfn]{Remark}

\begin{document}

\title{\vspace{-10pt}The number of spanning clusters of the uniform spanning tree in three dimensions}
\author{O.\ Angel\footnote{Department of Mathematics, University of British Columbia, Vancouver, BC, V6T 1Z2, Canada. Email: angel@math.ubc.ca.}, D.~A.~Croydon\footnote{Research Institute for Mathematical Sciences, Kyoto University, Kyoto 606-8502, Japan. Email: croydon@kurims.kyoto-u.ac.jp.}, S. Hernandez-Torres\footnote{Department of Mathematics, University of British Columbia, Vancouver, BC, V6T 1Z2, Canada. Email: saraiht@math.ubc.ca.}, and D.~Shiraishi\footnote{Department of Advanced Mathematical Sciences, Graduate School of Informatics,
Kyoto University, Kyoto 606-8501, Japan. Email: shiraishi@acs.i.kyoto-u.ac.jp.}}

\maketitle

\vspace{-30pt}

\begin{abstract}
Let ${\mathcal U}_{\delta}$ be the uniform spanning tree on $\delta \mathbb{Z}^{3}$. A spanning cluster of ${\mathcal U}_{\delta}$ is a connected component of the restriction of ${\mathcal U}_{\delta}$ to the unit cube $[0,1]^{3}$ that connects the left face $\{ 0 \} \times [0,1]^{2}$ to the right face $\{ 1 \} \times [0,1]^{2}$. In this note, we will prove that the number of the spanning clusters is tight as $\delta \to 0$, which resolves an open question raised by Benjamini in \cite{Itai}.
\end{abstract}

\vspace{-10pt}

\section{Introduction}

Given a finite connected graph $G = (V, E)$, a spanning tree $T$ of $G$ is a subgraph of $G$ that is a tree (i.e.\ is connected and contains no cycles) with vertex set $V$. A uniform spanning tree (UST) of $G$ is obtained by choosing a spanning tree of $G$ uniformly at random. This is an important model in probability and statistical physics, with beautiful connections to other subjects, such as electrical potential theory, loop-erased random walk and Schramm-Loewner evolution. See \cite{Lyons} for an introduction to various aspects of USTs.

Fix $\delta \in (0, 1)$ and $d\in\mathbb{N}$. In \cite{Pem} it was shown that, by taking the local limit of the uniform spanning trees on an exhaustive sequence of finite subgraphs of $\delta \mathbb{Z}^{d}$, it is possible to construct a random subgraph ${\mathcal U}_{\delta}$ of $\delta \mathbb{Z}^{d}$. Whilst the resulting graph ${\mathcal U}_{\delta}$ is almost-surely a forest consisting on an infinite number of disjoint components that are trees when $d \ge 5$, it is also the case that ${\mathcal U}_{\delta}$ is almost-surely a spanning tree of $\delta \mathbb{Z}^{d}$ with one topological end for $d \le 4$, see \cite{Pem}. In the latter low-dimensional case, $\mathcal{U}_\delta$ is commonly referred to as the UST on  $\delta \mathbb{Z}^{d}$.

In this note, we study a macroscopic scale property of ${\mathcal U}_{\delta}$, namely the number of its spanning clusters, as previously studied by Benjamini in \cite{Itai}. To be more precise, let us proceed to introduce some notation. Write
\begin{equation}\label{bdef}
\mathbb{B}  = [0, 1]^{d} = \left\{ (x_{1}, x_{2}, \cdots,  x_{d} ) \in \mathbb{R}^{d}\::\: 0 \le  x_{i} \le 1,\:i = 1, 2, \cdots , d \right\}
\end{equation}
for the unit hypercube in $\mathbb{R}^{d}$. Also, set
\begin{equation}\label{fdef}
F = \left\{ (x_{1}, x_{2}, \cdots,  x_{d} ) \in \mathbb{R}^{d}\::\:\ x_{1} = 0 \right\}
\end{equation}
and
\begin{equation}\label{gdef}
G = \left\{ (x_{1}, x_{2}, \cdots,  x_{d} ) \in \mathbb{R}^{d}\::\:\ x_{1} = 1 \right\}
\end{equation}
for the hyperplanes intersecting the `left' and `right' sides of the hypercube $\mathbb{B}$. Given a subgraph $U = (V, E) $ of $\delta \mathbb{Z}^{d}$, we write $U' = (V',E')$ for the restriction of $U$ to the cube $\mathbb{B} $, i.e.\ we set $V' = V \cap \mathbb{B}$ and $E' = \{ \{ x, y \} \in E :\: x, y \in V' \}$. A connected component of $U'$ is called a cluster of $U$. Moreover, following \cite{Itai}, a spanning cluster of $U$ is a cluster of $U$ containing vertices $x$ and $y$ such that $\text{dist}(x,F)<\delta$ and $\text{dist}(y,G)<\delta$, where $\text{dist}(z,A):=\inf_{w\in A}|z-w|$ is the Euclidean distance between a point $z\in \mathbb{R}^d$ and subset $A\subseteq\mathbb{R}^d$. That is, a cluster of $U$ is called spanning when it connects $F$ to $G$ (at the level of discretization being considered).

\begin{figure}[!t]
\begin{center}
\includegraphics[width=0.80\textwidth]{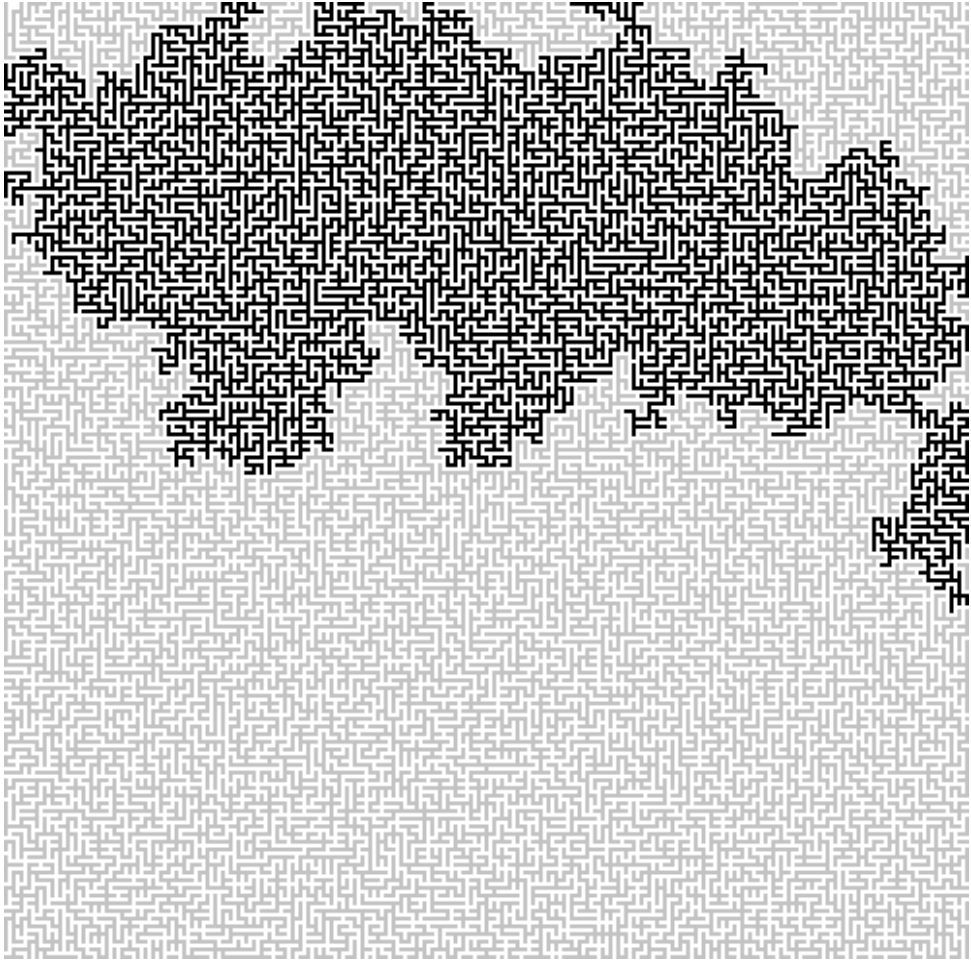}
\end{center}
\caption{Part of a UST in a two-dimensional box; the part shown is the central $115\times 115$ section of a UST on a $229\times229$ box. The single cluster spanning the two sides of the box is highlighted.}\label{spanningfig}
\end{figure}

Concerning the number of spanning clusters of ${\mathcal U}_{\delta}$, it was proved in \cite{Itai} that:
\begin{itemize}
\item for $d\geq 4$, the expected number of spanning clusters of ${\mathcal U}_{\delta}$ grows to infinity as $\delta \to 0$;
\item for $d=2$, the number of spanning clusters of ${\mathcal U}_{\delta}^{+}$ is tight as $\delta \to 0$, where ${\mathcal U}_{\delta}^{+}$ denotes the uniform spanning tree of the square $\mathbb{B} \cap \delta \mathbb{Z}^{2}$ when all the vertices on the right side of the square are identified to a single point (which is called the right wired uniform spanning tree in \cite{Itai}). Figure \ref{spanningfig} shows the spanning cluster of a realisation of (an approximation to) $\mathcal{U}_\delta$ on $\delta\mathbb{Z}^2$.
\end{itemize}
The case $d=3$ was left as an open question in \cite{Itai}. The main purpose of this note is to resolve it by showing the following theorem.

\begin{thm}\label{main}
Let $d=3$. It holds that the number of spanning clusters of ${\mathcal U}_{\delta}$ is tight as $\delta \to 0$.
\end{thm}

\begin{rem}
In the forthcoming article \cite{ACHS},
we establish a scaling limit for the three-dimensional UST in a version of the Gromov-Hausdorff topology, at least along the subsequence $\delta_n:=2^{-n}$. The corresponding two-dimensional result is also known (along an arbitrary sequence $\delta\rightarrow0$), see \cite{BCK} and \cite[Remark 1.2]{HS}. In both cases, we expect that the techniques used to prove such a scaling limit can be used to show that the number of spanning clusters of ${\mathcal U}_{\delta}$ actually converges in distribution. We plan to pursue this in a subsequent work that focusses on the topological properties of the three-dimensional UST.
\end{rem}

The organization of the remainder of the paper is as follows. In Section 2, we introduce some notation that will be used in the paper. The proof of Theorem \ref{main} is then given in Section 3.

\section{Notation}\label{notationsec}

In this section, we introduce the main notation needed for the proof of Theorem \ref{main}. We write $|\cdot|$ for the Euclidean norm on $\mathbb{R}^{3}$ and, as in the introduction, $\dist(\cdot, \cdot)$ for the Euclidean distance between a point and a subset of $\mathbb{R}^3$. Given $\delta \in (0,1)$, if $x\in\delta \mathbb{Z}^{3}$ and $r>0$, then we write
\[B(x,r) = \left\{ y \in \delta \mathbb{Z}^{3}:\:|x - y| < r\right\}\]
for the lattice ball of centre $x$ and radius $r$ (we will commonly omit dependence on $\delta$ for brevity). Let $\mathbb{B}$, $F$ and $G$ be defined as at \eqref{bdef}, \eqref{fdef} and \eqref{gdef} in the case $d=3$.

For $\delta\in(0,1)$, a sequence $\lambda = (\lambda (0), \lambda (1), \cdots , \lambda (m) )$ is said to be a path of length $m$ if $\lambda (i) \in \delta \mathbb{Z}^{3}$ and $|\lambda (i) - \lambda (i+1)| = \delta$ for every $i$. A path $\lambda$ is simple if $\lambda (i) \neq \lambda (j) $ for all $i \neq j$. For a path $\lambda = (\lambda (0), \lambda (1), \cdots , \lambda (m))$, we define its loop-erasure $\text{LE} (\lambda )$ as follows. Firstly, let
\[s_{0} = \max \left\{ j \le m:\:\lambda (j) = \lambda (0) \right\},\]
and for $i \ge 1$, set
\[s_{i} = \max \left\{ j \le m :\: \lambda (j) = \lambda ( s_{i-1} + 1 ) \right\}.\]
Moreover, write $n = \min \{ i :\: s_{i} = m  \}$. The loop-erasure of $\lambda$ is then given by
\[\text{LE} ( \lambda ) = \left( \lambda (s_{0} ), \lambda (s_{1} ), \cdots , \lambda (s_{n} ) \right).\]
We write $\text{LE} ( \lambda )  (k) = \lambda (s_{k})$ for each $0 \le k \le n$. Note that the vertices hit by $\text{LE} (\lambda )$ are a subset of those hit by $\lambda $, and that $\text{LE} (\lambda )$ is a simple path such that $\text{LE} (\lambda) (0) = \lambda (0)$ and $\text{LE} (\lambda) (n) = \lambda (m)$. Although the loop-erasure of $\lambda$ has so far only been defined in the case that $\lambda$ has a finite length, it is clear that we can define $\text{LE} (\lambda)$ similarly for an infinite path $\lambda$ if the set $\{ k \ge 0:\: \lambda (j) = \lambda (k) \}$ is finite for each $j \ge 0$. Additionally, when the path $\lambda$ is given by a simple random walk, we call $\text{LE}(\lambda) $ a loop-erased random walk (see \cite{Lawbound} for an introduction to loop-erased random walks).

Again given $\delta \in (0, 1)$, write ${\mathcal U}_{\delta}$ for the uniform spanning tree on $\delta \mathbb{Z}^{3}$. As noted in the introduction, this object was constructed in \cite{Pem}, and shown to be a tree with a single end, almost-surely. The graph ${\mathcal U}_{\delta}$ can be generated from loop-erased random walks by a procedure now referred to as Wilson's algorithm (after \cite{Wil}), which is described as follows.
\begin{itemize}
\item Let $(x_{i})_{i \ge 1}$ be an arbitrary, but fixed, ordering of $\delta \mathbb{Z}^{3}$.
\item Write $R^{x_{1}}$ for a simple random walk on $\delta \mathbb{Z}^{3} $ started at $x_{1}$. Let $\gamma_{x_{1}}  = \text{LE} ( R^{x_{1}}) $ be the loop-erasure of $R^{x_{1}}$ -- this is well-defined since $R^{x_{1}}$ is transient. Set ${\mathcal U}^{1} = \gamma_{x_{1}}$.
\item Given ${\mathcal U}^{i}$ for $i \ge 1$, let $R^{x_{i+1}}$ be a simple random walk (independent of ${\mathcal U}^{i}$) started at $x_{i + 1}$ and stopped on hitting ${\mathcal U}^{i}$. We let ${\mathcal U}^{i+1} = {\mathcal U}^{i} \cup \text{LE} (R^{x_{i+1}} ) $.
\end{itemize}
It is then the case that the output random tree $\cup_{i =1}^{\infty} {\mathcal U}^{i}$ has the same distribution as ${\mathcal U}_{\delta}$. In particular, the distribution of the output tree does not depend on the ordering of points $(x_{i})_{i \ge 1} $.

Similarly to above, for $z \in \delta \mathbb{Z}^{3}$, we will write $\gamma_{z}$ for the infinite simple path in ${\mathcal U}_{\delta}$ starting from $z$. Given a point $z \in \delta \mathbb{Z}^{3}$, it follows from the construction of ${\mathcal U}_{\delta}$ explained hitherto that the distribution of $\gamma_{z}$ coincides with that of $\text{LE} ( R^{z})$, where $R^{z}$ is a simple random walk on $\delta \mathbb{Z}^{3}$ started at $z$.

Furthermore, as we explained in the introduction, we will write ${\mathcal U}'_{\delta}$ for the restriction of ${\mathcal U}_{\delta}$ to the cube $\mathbb{B}$. A connected component of ${\mathcal U}'_{\delta}$ is called a cluster. Also, as we defined previously, a spanning cluster is a cluster connecting $F$ to $G$. We let $N_{\delta}$ be the number of spanning clusters of ${\mathcal U}_{\delta}$.

Finally, we will use $c$, $C$, $c_{0}$, etc.\ to denote universal positive constants which may change from line to line.

\section{Proof of the main result}

In this section, we will prove the following theorem, which incorporates Theorem \ref{main}.

\begin{thm} There exists a universal constant $C$ such that: for all $M < \infty$ and $\delta > 0$,
\begin{equation}\label{key}
\mathbf{P}\left( N_{\delta} \ge M \right) \le C M^{-1}.
\end{equation}
In particular, the laws of $(N_{\delta})_{\delta \in (0,1)}$ form a tight sequence of probability measures on $\mathbb{Z}_+$.
\end{thm}

\begin{rem}
In \cite{Aiz}, Aizenman proved that for critical percolation in two dimensions, the probability of seeing $M$ distinct spanning clusters is bounded above by $Ce^{- c M^{2}}$. We do not expect that the polynomial bound in \eqref{key} is sharp, but leave it as an open problem to determine the correct tail behaviour for number of spanning clusters of the UST in three dimensions, and, in particular, ascertain whether it also exhibits Gaussian decay.
\end{rem}

\begin{proof}
Let $\delta \in (0,1)$, and suppose $M \ge 1$ is such that $\delta < M^{-1}$. Define
\begin{align*}
A& = [-1, 2]^{3},\\
A'& = F \cap \mathbb{B},\\
A''& = \left\{ (x_{1}, x_{2}, x_{3} )\in\mathbb{B}:\:x_{1} = 2/3 \right\},\\
\mathbb{B}'& = \{ (x_{1}, x_{2}, x_{3} ) \in \mathbb{B} :\: x_{1} \le 2/3 \}.
\end{align*}
Moreover, let $(z_{i})_{i=1}^{L}$ be a sequence of points in  $A \cap \delta \mathbb{Z}^{3}$ such that
$A \subseteq \cup_{i=1}^{L} B (z_{i}, 1/M ) $ and $L \le  10^{5} M ^{3}$.

To construct ${\mathcal U}_{\delta}$, we first perform Wilson's algorithm for $(z_{i})_{i=1}^{L} $ (see Section \ref{notationsec}). Namely, we consider
\[{\mathcal U}^{1} := \bigcup_{i=1}^{L} \gamma_{z_{i} },\]
which is the subtree of ${\mathcal U}_{\delta}$ spanned by $(z_{i})_{i=1}^{L} $. (Recall that for $z \in \delta \mathbb{Z}^{3}$ we denote the infinite simple path in ${\mathcal U}_{\delta}$ starting from $z$ by $\gamma_{z}$.) The idea of the proof is then as follows. Crucially, each branch of ${\mathcal U}^{1}$ is a `hittable' set, in the sense that for a simple random walk $R$ whose starting point is close to ${\mathcal U}^{1}$, it is likely that $R$ hits ${\mathcal U}^{1}$ before moving far away. As a result, Wilson's algorithm guarantees that, with high probability, the spanning clusters of ${\mathcal U}_{\delta}$ correspond to those of ${\mathcal U}^{1}$ when $M$ is sufficiently large. So, the problem boils down to the tightness of the number of spanning clusters of ${\mathcal U}^{1}$, which is not difficult to prove.

To make the above argument rigorous, we introduce the following two ``good" events for ${\mathcal U}^{1}$:
\[H_{i} = H_{i} (\xi) := \left\{ \begin{array}{c}
\hspace{-5pt}\text{For any }x \in B(0,4)\cap\delta\mathbb{Z}^{3}\text{ with }\text{dist}(x,\gamma_{z_{i}})\le 1/M,\hspace{-5pt}\\
P^{x}_{R} \left( R[0, T ] \cap \gamma_{z_{i} } = \emptyset \right) \le M^{- \xi}
                                 \end{array}\right\},\]
\[I_{i} := \left\{ \begin{array}{c}
\hspace{-5pt}\text{The number of crossings of }\gamma_{z_{i} }\text{ between } A' \text{ and } A''\hspace{-5pt}\\
\text{ in } \mathbb{B}' \text{ is smaller than } M
                                 \end{array}\right\},\]
for $1 \le i \le L$, where
 \begin{itemize}
 \item $R$ is a simple random walk which is independent of $\gamma_{z_{i}}$, the law of which is denoted by $P^{x}_{R}$ when we assume $R(0) = x$;
 \item $T$ is the first time that $R$ exits $B ( x, 1/ \sqrt{M} )$;
 \item a crossing of $\gamma_{z_{i}} $ between $A'$ and $A''$ in $\mathbb{B}'$ is a connected component of the restriction of $\gamma_{z_{i}} $ to $\mathbb{B}'$ that connects $A'$ to $A''$.
\end{itemize}
Namely, the event $H_{i}$ guarantees that the branch $\gamma_{z_{i}}$ is a hittable set, and the event $I_{i}$ controls the number of crossings of $\gamma_{z_{i}}$.

Now, \cite[Theorem 3.1]{SS} ensures that there exist universal constants $\xi_{0}, C > 0$ such that
\[\mathbf{P} \left( \bigcap_{i=1}^{L} H_{i} (\xi_{0}) \right) \ge 1 - C M^{-10}.\]
Thus, with high probability (for ${\mathcal U}^{1}$), each branch of ${\mathcal U}^{1}$ is a hittable set.

The probability of the event $I_{i}$ is easy to estimate. Indeed, suppose that the event $I_{i}$ does not occur. This implies that the number of ``traversals'' of $S^{z_{i}}$ from $A'$ to $A''$ or vice versa must be bigger than $M$, where $S^{z_{i}}$ stands for a simple random walk starting from $z_{i}$. Notice that there exists a universal constant $c_{0}>0$ such that for any point $w \in A'$ (respectively $w \in A''$), the probability that $S^{w}$ hits $A''$ (respectively $A'$) is smaller than $1-c_{0}$ (see \cite[Proposition 1.5.10]{Lawb}, for example). Thus, the probability of the event $I_{i}$ is bounded below by $ 1- (1- c_{0} )^{M} =: 1 - e^{ - a M}$, where $a > 0$. Letting $b = a/2$ and taking sum over $1 \le i \le L$, we find that
\[\mathbf{P} \left( \bigcap_{i=1}^{L} I_{i}  \right) \ge 1 - C e^{ - b M}.\]

To put the above together, let
\[ J = \bigcap_{i=1}^{L} H_{i} (\xi_{0}) \cap I_{i}.\]
For $1 \le i \le L$, set ${\mathcal U}^{1}_{i} = \cup_{j=1}^{i} \gamma_{z_{j}}$ so that ${\mathcal U}^{1} = {\mathcal U}^{1}_{L}$. As above, by a spanning cluster of ${\mathcal U}^{1}_{i}$ between $A'$ and $A''$ in $\mathbb{B}'$ we mean a connected component of the restriction of ${\mathcal U}^{1}_{i}$ to $\mathbb{B}'$ which connects $A'$ to $A''$. We write $n_{i}$ for the number of spanning clusters of ${\mathcal U}^{1}_{i}$ between $A'$ and $A''$ in $\mathbb{B}'$. On the event $J$, we have that
\[n_{i} \le iM + i-1,\]
for all $1 \le i \le L$, since $n_{i+1} - n_{i}$ is at most $M+1$ for each $i \ge 1$. In particular, we see that the number of spanning clusters of ${\mathcal U}^{1}$ between $A'$ and $A''$ in $\mathbb{B}'$ is bounded above by $L (M+1)$, which is comparable to $M^{4}$.

We next consider a sequence of subsets of $A$ as follows. Let $a^{\ast} > 0$ be the positive constant such that
\begin{equation}\label{const}
a^{\ast} \sum_{k=1}^{\infty} k^{-2} = 10^{-1}.
\end{equation}
Set $\eta_{1} = 0$, and $\eta_{k} = a^{\ast} \sum_{j=1}^{k-1} j^{-2}$ for $k \ge 2$. Finally,  for $k \ge 1$, let
\[A_{k} = [-1+ \eta_{k} , 2- \eta_{k} ]^{3}.\]
Notice that $A_{k+1} \subseteq A_{k}$ and $[-1/2, 3/2]^{3} \subseteq A_{k}$ for all $k \ge 1$, and moreover $\text{dist}(\partial A_{k},\partial A_{k+1}) = a^{\ast} k^{-2} $. We further introduce sequences $(z_{i}^{k})_{i =1}^{L_{k}}$ consisting of points in $A_{k} \cap \delta \mathbb{Z}^{3}$ such that
\[A_{k} \subseteq \bigcup_{i=1}^{L_{k}} B \left( z_{i}^{k}, \delta_{k} \right),\]
and $L_{k} \le 10^{5}  \delta_{k}^{-3}$, where $\delta_{k} := M^{-1} 2^{-(k-1)}$. Note that we may assume that $L_{1} = L$ and $(z_{i}^{1})_{i =1}^{L_{1}} = ( z_{i})_{i=1}^{L} $.

We set
\[H_{i}^{k} = H_{i}^{k} (\xi) := \left\{
\begin{array}{c}
 \hspace{-5pt}\text{For any } x \in B (0, 4) \cap \delta \mathbb{Z}^{3} \text{ with } \text{dist} \left( x, \gamma_{z_{i}^{k} } \right) \le \delta_{k},\hspace{-5pt}\\
   P^{x}_{R} \left( R[0, T^{k} ] \cap \gamma_{z_{i}^{k} } = \emptyset \right) \le \delta_{k}^{ \xi}
\end{array}\right\},\]
where $R$ is a simple random walk that is independent of $\gamma_{z_{i}^{k}}$, with law denoted by $P^{x}_{R}$ when we assume $R(0) =x$, and $T^{k}$ is the first time that $R$ exits $B(x,\sqrt{\delta_{k}})$. By \cite[Theorem 3.1]{SS} again, there exist universal constants $\xi_{1}, C > 0$ (which do not depend on $k$) such that
\[\mathbf{P}\left( \bigcap_{i=1}^{L_{k}} H_{i}^{k} (\xi_{1}) \right) \ge 1 - C \delta_{k}^{10},\]
for all $k =1, 2, \cdots , k_{0}$, where $k_{0}$ is the smallest integer $k$ such that $\delta_{k} < \delta$. Thus if we write
\[H^{k} =  \bigcap_{i= 1}^{L_{k}} H_{i}^{k} (\xi_{1})\]
and
\[J' = J \cap  \bigcap_{k=1}^{k_{0}} H^{k},\]
we have
\[\mathbf{P} (J') \ge 1 - C M^{-10}.\]

Given the above setup, we perform Wilson's algorithm as follows:
\begin{itemize}
\item recall that ${\mathcal U}^{1}$ is the tree spanned by $(z_{i}^{1})_{i =1}^{L_{1}} = (z_{i})_{i=1}^{L}$;
\item next perform Wilson's algorithm for $(z_{i}^{2} )_{i =1}^{L_{2}}$ --  for each $z_{i}^{2}$, run a simple random walk $R^{z_{i}^{2}}$ from $z_i^2$ until it hits the part of the tree that has already been constructed, and adding its loop-erasure as a new branch -- the output tree is denoted by ${\mathcal U}^{2}$;
\item repeat the previous step for $(z_{i}^{k})_{i =1}^{L_{k}}$ to construct ${\mathcal U}^{k}$ for $k= 1,2, \cdots , k_{0}$.
\end{itemize}

Now, condition ${\mathcal U}^{1}$ on the event $J$ above. We will show that, with high (conditional) probability, every new branch in ${\mathcal U}^{2} \setminus {\mathcal U}^{1}$ has diameter smaller than $M^{-1/4}$. To this end, for $1 \le i \le L_{2}$, we write $d^{2}_{i}$ for the Euclidean diameter of the path from $z_{i}^{2}$ to ${\mathcal U}^{1}$ in ${\mathcal U}^{2}$, and define the event $W^{2}_{i}$ by setting
\[W^{2}_{i} = \left\{ d^{2}_{i} \ge M^{-1/4} \right\}.\]
Suppose that the event $W^{2}_{i}$ occurs. By Wilson's algorithm, the simple random walk $R^{z_{i}^{2}}$ must not hit ${\mathcal U}^{1}$ until it exits $B (z_{i}^{2}, M^{-1/4})$. Since $\text{dist} (z_{i}^{2}, \partial A ) \ge a^{\ast}$ (for the constant $a^{\ast}$ defined at \eqref{const}), it holds that $B (z_{i}^{2}, M^{-1/4}) \subseteq A$. With this in mind, we set $u_{0} = 0$, and
\[u_{m} =  \inf \left\{ j \ge u_{m-1}:\: \left| R^{z_{i}^{2}} (j) - R^{z_{i}^{2}} (u_{m-1}) \right| \ge M^{-1/2} \right\}\]
for $m \ge 1$. We then have that
\[R^{z_{i}^{2}} [ u_{m-1} , u_{m } ] \cap {\mathcal U}^{1} = \emptyset\]
for all $1 \le m \le M^{1/4}$. Since $A \subseteq \cup_{i=1}^{L} B (z_{i}, 1/M)$, it follows that for each $1 \le m \le M^{1/4}$, there exists a $z_{i}$ such that $ R^{z_{i}^{2}} ( u_{m-1} ) \in B (  z_{i} , 1/M )$. Thus the event $H_{i}(\xi_{0})$ guarantees that
\[\mathbf{P} \left( R^{z_{i}^{2}} [ u_{m-1} , u_{m } ] \cap {\mathcal U}^{1} = \emptyset \text{ for all } 1 \le m \le M^{1/4} \right) \le M^{- \xi_{0} M^{1/4}}.\]
Consequently, the conditional probability of $\cup_{i=1}^{L_{2}} W^{2}_{i}$ is bounded above by $L_{2} M^{- \xi_{0} M^{1/4}}$, which is smaller than $Ce^{- c M^{1/4}}$. This implies that,  with probability at least $ 1- Ce^{- c M^{1/4}}$, every new branch in ${\mathcal U}^{2} \setminus {\mathcal U}^{1}$ has diameter smaller than $M^{-1/4}$. Notice that once each new branch has such a small diameter, the event $J$ guarantees that the number of spanning clusters of ${\mathcal U}^{2}$ between $A'$ and $A''_{2}$ in $\mathbb{B}$ is bounded above by $L (M+1) \le 10^{6} M^{4} $, where $A''_{2}$ is defined by setting
\[A''_{2} = \mathbb{B} \cap \left\{(x_1,x_2,x_3)\in\mathbb{R}^3:\:x_{1} = 2/3 + M^{-1/4}\right\}.\]

Essentially the same argument is valid for ${\mathcal U}^{k}$. Indeed, conditioning ${\mathcal U}^{k}$ on the good event $J \cap \cap_{l=1}^{k} H^{l}$ as above, it holds that, with probability at least $1 - Ce^{- c \delta_{k}^{- 1/ 4}}$ every new branch in ${\mathcal U}^{k+1} \setminus {\mathcal U}^{k}$ has diameter smaller than $\delta_{k}^{1/4}$. Notice that $\sum_{k} \delta_{k}^{1/4} \le 10 M^{-1/4} < 10^{-2}$ when $M$ is large. Therefore, with probability at least $1 - C M^{-10}$, the number of spanning clusters of ${\mathcal U}^{k_{0}}$ between $A'$ and $A'''$ in $\mathbb{B}$ is bounded above by $L (M+1) \le C M^{4} $ for some universal constant $C$, where $A'''$ is defined by setting
\[A''' = \mathbb{B} \cap \left\{(x_1,x_2,x_3)\in\mathbb{R}^3:\: x_{1} = 3/4 \right\}.\]

Finally, we perform Wilson's algorithm for all of the remaining points in $\delta \mathbb{Z}^{3}$ to construct ${\mathcal U}_{\delta}$. Since the mesh size of $( z_{i}^{k_{0}} )_{i =1}^{L_{k_{0}}}$ is smaller than $\delta$, it follows that the restriction of ${\mathcal U}_{\delta}$ to $\mathbb{B}$ coincides with that of ${\mathcal U}^{k_{0}}$. Thus we conclude that there exists a universal constant $C$ such that: for all $M < \infty$ and $\delta \in (0, M^{-1})$,
\[\mathbf{P} \left( N_{\delta} \ge M \right) \le C M^{-2}.\]
For the case that $\delta > M^{-1}$, it is clear that $N_{\delta} < 100 M^{2}$. Combining these two bounds, we readily obtain the bound at \eqref{key}.
\end{proof}

\section*{Acknowledgements}

DC would like to acknowledge the support of a JSPS Grant-in-Aid for Research Activity Start-up, 18H05832 and a JSPS Grant-in-Aid for Scientific Research (C), 19K03540. SHT is supported by a fellowship from the Mexican National Council for Science and Technology (CONACYT). DS is supported by a JSPS Grant-in-Aid for Early-Career Scientists, 18K13425.

\providecommand{\bysame}{\leavevmode\hbox to3em{\hrulefill}\thinspace}
\providecommand{\MR}{\relax\ifhmode\unskip\space\fi MR }
\providecommand{\MRhref}[2]{%
  \href{http://www.ams.org/mathscinet-getitem?mr=#1}{#2}
}
\providecommand{\href}[2]{#2}

\end{document}